\documentclass{jT}
\usepackage{amssymb,amsmath,latexsym,amsthm}

\def\Pari{\hbox{\textsc{PARI/GP}}}
\def\Magma{\hbox{\textsc{Magma}}}
\def\Macsyma{\hbox{\textsc{Macsyma}}}
\def\Q{\hbox{\textbf{Q}}}

\begin{document}
\title{A note on integral points on elliptic curves}
\author[Mark {\sc Watkins}]{{\sc Mark} WATKINS}
\address{Mark {\sc Watkins}\\
Department of Mathematics\\
University Walk\\
University of Bristol\\
Bristol, BS8 1TW\\
England}
\email{watkins@maths.usyd.edu.au}
\maketitle

\begin{abstr}
We investigate a problem considered by Zagier and Elkies,
of finding large integral points on elliptic curves.
By writing down a generic polynomial solution and equating coefficients,
we are led to suspect four extremal cases
that still might have nondegenerate solutions.
Each of these cases gives rise to a polynomial system of equations,
the first being solved by Elkies in 1988 using the resultant methods
of~\Macsyma, with there being a unique rational nondegenerate solution.
For the second case we found that resultants and/or Gr\"obner bases
were not very efficacious.
Instead, at the suggestion of Elkies, we used multidimensional
$p$-adic Newton iteration, and were able to find a nondegenerate
solution, albeit over a quartic number field. Due to our methodology,
we do not have much hope of proving that there are no other solutions.
For the third case we found a solution in a nonic number field,
but we were unable to make much progress with the fourth case.
We make a few concluding comments and include an appendix from Elkies
regarding his calculations and correspondence with Zagier.
\end{abstr}

\begin{resume}
\`A la suite de Zagier et Elkies, nous recherchons de grands points entiers
sur des courbes elliptiques. En \'ecrivant une solution polynomiale
g\'en\'erique et en \'egalisant des coefficients, nous obtenons quatre cas
extr\'emaux susceptibles d'avoir des solutions non d\'eg\'en\'er\'ees.
Chacun de ces cas conduit \`a un syst\`eme d'\'equations polynomiales,
le premier \'etant r\'esolu par Elkies en 1988 en utilisant les
r\'esultants de Macsyma; il admet une unique solution rationnelle non
d\'eg\'en\'er\'ee.
Pour le deuxi\`eme cas nous avons constat\'e que les r\'esultants ou les
bases de Gr\"obner sont peu efficaces.  Suivant une suggestion d'Elkies,
nous avons alors utilis\'e une it\'eration de Newton $p$-adique
multidimensionnelle et d\'ecouvert une solution non d\'eg\'en\'er\'ee,
quoique sur un corps de nombres quartique. En raison de notre
m\'ethodologie, nous avons peu d'espoir de montrer qu'il n'y a aucune
autre solution. Pour le troisi\`eme cas nous avons trouv\'e une solution sur
un corps de degr\'e 9, mais n'avons pu traiter le quatri\`eme cas. Nous
concluons par quelques commentaires et une annexe d'Elkies concernant
ses calculs et sa correspondance avec Zagier.
\end{resume}

\section{Introduction}
Let $E$ be an elliptic curve given by the model~$y^2=x^3+Ax+B$,
and suppose that $(X,Y)$ is an integral point on this model.
How large can $X$ be in terms of $|A|$ and~$|B|$?
One measure of the impressiveness of the size
of an integral point is given by the quotient
$\rho=\log(X)/\log\bigl(\max(|A|^{1/2},|B|^{1/3})\bigr)$,
which, as Zagier~\cite{zagier} indicates, can be interpreted as
saying that $X$ is of the order of magnitude of the $\rho$th power of the
roots of the cubic polynomial~$x^3+Ax+B$.

Lang \cite{lang} makes the conjecture that $\rho$ is bounded,
and notes (see \cite{zagier}) that he and Stark worked out that
generically $\rho\le 10+o(1)$ via probablistic heuristics,
though a construction of Stark indicated that in similar
situations there might be finitely many exceptional parametric families
with larger~$\rho$. Vojta \cite{vojta} has related this conjecture
to his more general Diophantine theory, where again these exceptional
families cannot be eliminated.
In 1987, Zagier \cite{zagier} gave a construction that gives
infinitely many curves with $\rho\ge 9-o(1)$, and listed some
impressive examples from numerical calculations of Odlyzko.

In a letter to Zagier in 1988, Elkies constructed infinitely
many examples that satisfy~$\rho\ge 12-o(1)$.
His construction is polynomial-based, and reduces to solving a
system of polynomial equations formed from equating coefficients.
There are exactly four choices of parameters that both yield $\rho=12-o(1)$
and for which there is a reasonable hope that a solution might exist.
The first of these was the case worked out by Elkies.
This already led to a system of 4 polynomial equations in 4 variables,
which Elkies notes took a longish session
of \Macsyma~\cite{macsyma-ref} to solve.
The second choice of parameters immediately
(via linear substitution) leads to a system of 6 equations and unknowns;
even though computers have gained much in speed over the last 15 years,
the resulting system is still too difficult to solve
via Gr\"obner bases or resultants.
We eliminated one variable from the system via another linear substitution
(though this creates denominators), and then another via a resultant step.
This gives us a rather complicated system of four equations and unknowns;
the degrees of the polynomials were sufficiently large that, again,
Gr\"obner bases and resultants were not of much use.
We then proceeded to try to find solutions via a multidimensional
\hbox{$p$-adic} iterative Newton method. We found one such solution
over a quartic number field; it is an inherent problem with this method
that we have little hope of proving that we have found all the solutions.
With the third choice of parameters, we found a solution in a nonic
number field, and with the fourth case we made little progress.

As an appendix, we include some calculations of Elkies regarding the
first case, and his 1988 letter to Zagier.

\section{Families of Pell type}
First we review the construction of Elkies. Consider the equation
\begin{equation}
X(t)^3+A(t)X(t)+B(t)=Q(t)Y(t)^2\label{eqn:EPZ}
\end{equation}
where $A,B,Q,X,Y$ are polynomials in~$t$ with $Q$ quadratic.
Given a rational polynomial solution to this equation, via scaling
we can make all the polynomials integral. The theory of the Pell
equation implies that if the quadratic polynomial $Q(t)$
is a square for one integral $t$-value,
then it is square for infinitely many integral~$t$,
and thus we get infinitely many curves $y^2=x^3+A(t)x+B(t)$
with integral points~$\bigl(X(t),Y(t)\sqrt{Q(t)}\bigr)$.

Let $a,b,q,x,y$ be the degrees of these polynomials respectively.
We wish for $\rho=x/\max(a/2,b/3)$ to be as large as possible.
If we do a parameter count, we get that there are $(a+b+q+x+y)+5$
coefficients of our polynomials. The total degree of our equation
is $3x=q+2y$, so we get $3x+1$ equations.
When $3x+1\le a+b+q+x+y+5$, we might expect there to be a solution.
However, we first need to remove the effect of the action of the
group ${\rm PGL}_2(\Q)$ on our choice of coefficents.

Letting $l(P)$ be the leading coefficient of a polynomial~$P$,
we first scale $t$ by $l(X)/l(Y)$ and then multiply through $l(Y)^x/l(X)^y$,
so as to make $X,Y,Q$ all monic.
Then we translate so as to eliminate the $t^{y-1}$ term in $Y$.
Then we effect $t\rightarrow 1/t$ and multiply $(X,Y,Q,A,B)$
by $(t^x,t^y,t^2,t^{2x},t^{3x})$, and then scale so as to make\footnote{
 We could alternatively equate two coefficients, or set the linear
 coefficient of $Q$ equal to~1; we found that fixing the linear coefficient
 of $X$ to be~1 was best amongst the various choices.
 In this scaling, we assume the coefficient is nonzero;
 the alternative case can be handled separately.}
the $t$-coefficient of $X$ be equal to~1.
Finally we undo the $t\rightarrow 1/t$ transformation in the same manner.
So we are left with $(a+1)+(b+1)+q+(x-1)+(y-1)$ coefficients,
while we also lose one condition, namely that the leading coefficients match.
Thus we want to have $a+b+q+x+y\ge 3x$ with $\rho=x/\max(a/2,b/3)$
as large as possible, and this turns out to be~12.
We get 4 different possibilities, namely
$(a,b,q,x,y)=(0,1,2,4,5),(1,1,2,6,8),(1,2,2,8,11),(2,3,2,12,17)$.
For instance, for the first case we have the polynomials
\begin{align*}
X(t)&=t^4+t^3+x_2t^2+x_1t+x_0,\quad Y(t)=t^5+y_3t^3+y_2t^2+y_1t+y_0,\\
&Q(t)=t^2+q_1t+q_0,\quad A(t)=a_0,\quad B(t)=b_1t+b_0,
\end{align*}
and equating the $t^0$-$t^{11}$ coefficients gives us 12 equations in
these 12 unknowns. Fortunately, simple linear substitutions easily reduce
this to 4 equations and unknowns; we give one such reduced set, in order
to indicate the complexity of the equations.
\begin{align*}
&\scriptstyle
{12x_0x_2-12x_0q_0+60x_0+6x_1^2-24x_1x_2+48x_1q_0-156x_1-
x_2^3-3x_2^2q_0+27x_2^2+9x_2q_0^2-174x_2q_0+}\\
&\>\scriptstyle{+417x_2-5q_0^3+171q_0^2 -939q_0+1339=0,}
\end{align*}
\begin{align*}
&\scriptstyle{
4x_0x_1+4x_0x_2+4x_0q_0+8x_0+2x_1^2-x_1x_2^2-2x_1x_2q_0-6x_1x_2+
3x_1q_0^2-10x_1q_0-17x_1+2x_2^2q_0+5x_2^2-}\\
&\>\scriptstyle{-12x_2q_0^2+26x_2q_0+38x_2+10q_0^3-71q_0^2+80q_0+83=0,}
\end{align*}
\begin{align*}
&\scriptstyle{
120x_0x_1x_2-72x_0x_1q_0+312x_0x_1-60x_0x_2^2+216x_0x_2q_0-576x_0x_2-
60x_0q_0^2+336x_0q_0-516x_0+}\\
&\>\scriptstyle{+32x_1^3-168x_1^2x_2+288x_1^2q_0-936x_1^2-
18x_1x_2^3-54x_1x_2^2q_0+342x_1x_2^2+
114x_1x_2q_0^2-1836x_1x_2q_0+}\\
&\>\>\scriptstyle{+4146x_1x_2-42x_1q_0^3+1302x_1q_0^2-6870x_1q_0
+9642x_1+9x_2^4+72x_2^3q_0-234x_2^3-342x_2^2q_0^2+2658x_2^2q_0-}\\
&\>\>\scriptstyle{-4488x_2^2+336x_2q_0^3-4518x_2q_0^2+17004x_2q_0
-19446x_2-75q_0^4 +1486q_0^3-9036q_0^2+22098q_0=19041,}
\end{align*}
\begin{align*}
&\scriptstyle{
48x_0^2x_2-48x_0^2q_0+240x_0^2+64x_0x_1^2-128x_0x_1x_2+288x_0x_1q_0
-768x_0x_1-40x_0x_2^3+24x_0x_2^2q_0-}\\
&\>\scriptstyle{
-72x_0x_2^2+40x_0x_2q_0^2-816x_0x_2q_0+1864x_0x_2-
24x_0q_0^3+792x_0q_0^2-4232x_0q_0+5928x_0-28x_1^2x_2^2-}\\
&\>\scriptstyle{-24x_1^2x_2q_0+24x_1^2x_2
+36x_1^2q_0^2-312x_1^2q_0 +    
420x_1^2+84x_1x_2^3-84x_1x_2^2q_0+388x_1x_2^2-
244x_1x_2q_0^2+}\\
&\>\scriptstyle{+1480x_1x_2q_0-1876x_1x_2+180x_1q_0^3-2028x_1q_0^2+
6428x_1q_0-6180x_1+3x_2^5+9x_2^4q_0-84x_2^4-26x_2^3q_0^2+}\\
&\>\>\scriptstyle{+480x_2^3q_0-1186x_2^3+10x_2^2q_0^3-224x_2^2q_0^2+
1310x_2^2q_0-2140x_2^2+7x_2q_0^4-464x_2q_0^3+4210x_2q_0^2-}\\
&\>\>\scriptstyle{-12472x_2q_0+11807x_2-3q_0^5 +228q_0^4-2942q_0^3+
14284q_0^2-29791q_0+22560=0.}
\end{align*}

If we are willing to accept variables in denominators, we can go one
step more and eliminate $x_0$ from one of the first three equations.
A system like this was solved by Elkies in 1988 using \Macsyma\ which
uses resultants; solving it is almost instantaneous\footnote{
 That is, provided one deals with the multivariate polynomial rings
 properly and works over the rationals/integers at the desired times.}
with \Magma\ today, using either Gr\"obner bases or resultants.
We get an isolated solution and also two (extraneous) positive-dimensional
solution varieties (which correspond to points
on the singular plane cubic curve):
\begin{align*}
(x_0,x_1,x_2,q_0)=&
\textstyle
{\bigl({1\over 192}[16u^2-200u-239],{1\over 8}[4u-1],u,{9\over 4}\bigr),}\\
&\qquad\qquad\qquad\qquad(u,-2v+3,v-5,v),
\textstyle{({311\over 64},{61\over 8},{9\over 2},{11\over 4})}.
\end{align*}
From the isolated point, via back-substitution we get
$$(y_0,y_1,y_2,y_3,q_1,a_0,b_0,b_1)=
\textstyle{({715\over 64},{165\over 16},{77\over 16},{55\over 8},
3,{216513\over 4096},-{3720087\over 131072},{531441\over 8192}).}$$

To derive the solution in the form given by Elkies, we first want to
eliminate denominators, and we also wish to minimise the value of~$A$
that occurs at the end (that is, get rid of spurious powers of $2$ and~$3$).
This can be done by replacing $t$ by $1-9t/2$
and then multiplying $(X,Y,Q,A,B)$ by $(s,-4s/3,9s/16,s^2,s^3)$
where $s=128/81$. This gives us
\begin{align*}
X(t)&=6(108t^4-120t^3+72t^2-28t+5),\\
Y(t)&=72(54t^5-60t^4+45t^3-21t^2+6t-1),\\
Q(t)&=2(9t^2-10t+3),\> A(t)=132,\> B(t)=-144(8t-1).
\end{align*}
Note that $Q(1)=2^2$, so that there are infinitely many integral
values of~$t$ for which $Q(t)$ is square. As noted by Elkies,
we have that $X(t)\sim B(t)^4/2^{25}3^4$, so that small values
of~$t$ do not give very impressive values of~$\rho$.

\def\EPZII{\hbox{${\rm EPZ}_{\rm II}$}}
\def\RII{\hbox{${\rm R}_{\rm II}$}}
\subsection{The second case}
We next consider the second case \EPZII\ of the Elkies-Pell-Zagier
equation~\eqref{eqn:EPZ}, where $(a,b,q,x,y)=(1,1,2,6,8)$.
After making rational transformations, we are left with 18 equations
in 18 unknowns, which reduce to 6 upon making linear substitutions.
We can reduce to 5 via allowing denominators,\footnote{
 This linear substitution is probably most efficiently done via resultants,
 as else the denominators will cause problems for some
 computer algebra systems.}
and then eliminate one more variable via resultants,
but at this point, we are left with
equations with too large of degrees for resultants or Gr\"obner bases
to be of much use. Parts of two of the four equations appear below
(the whole input file is about 500 kilobytes)
$$2101324894157987694q_0^{14}+
107129273851487767680x_2^2x_3^2x_4^2q_0^5+\cdots=0,$$
\begin{align*}
&32970900880723713844451225823q_0^{22}-\\
&\qquad\qquad-34328441295817679913295188031488
x_2^2x_3^2x_4^7q_0^6+\cdots=0.
\end{align*}
We denote this reduced system of equations by~\RII.

It was suggested to us by Elkies that it might be
possible to find a solution via multidimensional $p$-adic Newton
iteration.\footnote{
 This technique appears in~\cite{elkies},
 while J.~Wetherell tells us that he has used it to find torsion
 points on abelian varieties. In \cite{elkies}, the lifting step
 was done via computing derivatives numerically, while we chose
 to compute them symbolically. Uses of this technique in situations
 close to those that occur here will be described in~\cite{antsVII}.}
In general, this method is most useful when we are searching for
zero-dimensional solution varieties in a small number of variables.
Writing $\vec f$ as our system of equations, we take a $p$-adic
approximate solution $\vec s$ and replace
it by $\vec s-J(\vec s)^{-1}\vec f(\vec s)$,
where $J(\vec s)$ is the Jacobian matrix of partial derivatives
for our system evaluated at~$\vec s$.
Since the convergence is quadratic, it is not
difficult to get $p$-adic solutions to high precision.
From each liftable local solution mod~$p$
we thus obtain a solution modulo a large power of~$p$,
and then use standard lattice reduction techniques \cite[\S 2.7.2]{Cohen}
to try to recognise it as a rational or algebraic number.

First we tried the primes $p=2,3$, but we found no useful mod~$p$ solutions;
all the local solutions had a noninvertible Jacobian matrix.\footnote{
 Many of them had a Jacobian equal to the zero matrix, and these we
 expect to come from positive-dimensional solution varieties.}
Furthermore, since a solution to \RII\ might very well have coordinates
whose denominators have powers of 2 and~3, not finding a solution was not
too surprising.
With $p=5$ we again found some (probable) positive-dimensional families
and three other solutions, of which two had an invertible Jacobian modulo~5.
However, these solutions to \RII\ failed to survive the undoing of the
resultant step, and thus do not actually correspond to a solution to~\EPZII.
We found the same occurrence for $p=7,11,13$ --- there were various $\Q_p$
solutions to our reduced system,
but these did not lift back to original system.

With $p=17$ our luck was better, as here we found a solution in the
dihedral quartic number field $K$ defined by~$z^4-2z^3-4z^2+5z-2$,
whose discriminant is $-3^2 11^3$.
Letting $\theta$ be a root of this polynomial,
the raw form of our solution is
$$x_2=\textstyle{1\over 2430000}
(9069984\theta^3+66428384\theta^2+19934816\theta-283298787),$$
$$x_3=\textstyle{1\over 6750}
(20240\theta^3+70576\theta^2-121616\theta-441839),$$
$$x_4=\textstyle{1\over 900}
(-5808\theta^3-7568\theta^2+33968\theta+23959),$$
$$q_0=\textstyle{1\over 2700}
(2576\theta^3+3760\theta^2-8720\theta+10971).$$
After undoing the resultant step the rest is but substitution
and we readily get a solution to~\EPZII,
albeit, in a quartic number field. Note that our prime 17 is the
smallest odd unramified prime which has a degree 1 factor in~$K$; due to our
method of division of labour we actually first found the solution mod~29.
Since we do not know $K$ ahead of time, we have little choice
but to try all small primes.

We next introduce some notation before stating our result;
we have infinitely many Pell equations from which to choose,
and so only present the simplest one that we were able to obtain.
Let $$p_2=\theta,\quad q_2=\theta-1,\quad r_2=\theta^2-\theta-5,
\quad\text{and}\quad p_3=2\theta^2-2\theta+1 $$
be the primes above~2 and the ramified prime above~3,
and $$\eta_1=\theta^3+\theta^2-2\theta+1 \quad\text{and}\quad
\eta_2=\theta^3-3\theta+1$$ be fundamental units, so that
we have $p_2q_2r_2=2$ and $p_3^2=3\eta_1^2\eta_2^{-1}$.
Let $\beta=2\theta^3+2\theta^2-6\theta-3$ (this is of norm~3271),
and with
$$Q(t)=c_2t^2+c_1t+c_0=
3p_2^7q_2\beta\eta_1^2\eta_2^{-1} t^2+
2q_2^3\eta_1^2\beta(\theta^3-\theta^2+11)t+
q_2^2\beta^2\eta_2^2$$ we have
\begin{align*}
X(t)&=2^4 3^4p_2^5q_2^7\beta\eta_1^8\eta_2^{-4} t^6+
2^3 3^4q_2^5r_2\beta\eta_1^9\eta_2^{-4}
(17\theta^3+2\theta^2-71\theta+33)t^5+\\
&\quad +2^2 3^3q_2\beta\eta_1^8\eta_2^{-3}
(1463\theta^3-2436\theta^2-2667\theta+1903)t^4+\\
&\quad\quad +24q_2\beta\eta_1^6\eta_2^{-2}
(25901\theta^3+32060\theta^2-52457\theta+15455)t^3+\\
&\quad\quad\quad +12q_2^2p_3\beta\eta_1^3
(40374\theta^3+47422\theta^2-61976\theta+37707)t^2+\\
&\quad\quad\quad\quad +2q_2^2r_2^2p_3\beta\eta_1^3\eta_2
(7081\theta^3-854\theta^2+90791\theta-23035)t+\\
&\quad\quad\quad\quad\quad +q_2\beta\eta_1\eta_2^2
(190035\theta^3+199008\theta^2-174189\theta+50449),
\end{align*}
\begin{align*}
A(t)&=
-12q_2^4r_2p_3\beta^2\eta_1\eta_2^{-3}(\theta^3-\theta+1)t-\\
&\quad -q_2^2p_3\beta^2\eta_1^{-1}\eta_2^{-2}
(\theta^3-\theta+1)(9\theta^3-2\theta^2+5\theta+9),
\end{align*}
and
\begin{align*}
B(t)&=-6q_2^7r_2^2\beta^3\eta_1\eta_2^{-4}(2\theta-1)^4t-\\
&\quad -q_2^4r_2\beta^3\eta_1^{-1}\eta_2^{-3}(2\theta-1)^4
(4\theta^3+18\theta^2-16\theta+1).
\end{align*}
With the above definition of $c_2=3p_2^7q_2\beta\eta_1^2\eta_2^{-1}$,
we have that
$$f_1=p_2^{-1}q_2p_3^{-1}\eta_1^{-3}\eta_2
(\theta^3+2\theta^2-\theta+1)\sqrt{c_2}+
r_2\eta_1^{-1}\eta_2(3\theta^3-19\theta^2+20\theta-5),$$
$$f_2=2^2p_2^4p_3^{-1}\eta_1^{-1}\eta_2^{-1}\sqrt{c_2}+
\eta_1\eta_2^{-1}(19\theta^3-51\theta^2+38\theta-5),$$
and
$$f_3=p_2q_2^2\eta_1^{-4}\eta_2^3(6\theta^2-2\theta+1)\sqrt{c_2}+
r_2\eta_1\eta_2(19\theta^3-14\theta^2-71\theta-41)$$
are units of relative norm~1 in $K\bigl(\sqrt{c_2}\bigr)$.
Again from the above definitions we have
$\sqrt{c_0}=\pm q_2\beta\eta_2=\pm(49\theta^3+41\theta^2-77\theta+33),$
and so we solve the Pell equation and obtain square values of $Q(t)$ by taking
$$t=2\sqrt{c_0}uv+v^2c_1\quad
\text{where}\quad f_1^if_2^jf_3^k=u+v\sqrt {c_2}$$
for integers $i,j,k$.
We can make $t$ integral via various congruence restrictions
on $(i,j,k)$; however, note that $p_2$ divides all but the constant
coefficients of our polynomials (including~$Y$),
and so we still get integral solutions
to~\EPZII\ even when $p_2$ exactly divides the denominator of~$t$.
Similarly, the only nonconstant coefficient that~$p_3$ fails to divide
is the linear coefficient of~$Q$; since $p_3^3$ divides $Y(t)$,
this nuisance evaporates when we consider solutions to~\EPZII.
As $t\rightarrow\infty$ the norm of the ratio $X(t)/A(t)^6$
tends to $1/2^{56}3^{20}17^6 3271^{11}$; we do not know if this
is as large as possible.

We are fairly certain that there are no more nondegenerate
algebraic solutions to~\EPZII, but we have no proof of this.
For the small primes, we have identified every local solution
to~\RII\ that has invertible Jacobian as algebraic.
In addition to the above quartic solution, there are five such
solutions\footnote{
 Since the local images of these solutions to~\RII\ failed to survive
 the undoing of the resultant step modulo~$p$,
 this determination of their algebraicity is unnecessary
 as evidence toward our claim that \EPZII\ has no more solutions,
 but might be interesting in that it shows the splitting of a
 large-dimensional algebra into many smaller fields.}
having degrees 13, 17, 19, 22, and~22, with each having maximal Galois group.

\def\EPZIII{\hbox{${\rm EPZ}_{\rm III}$}}
\def\EPZIV{\hbox{${\rm EPZ}_{\rm IV}$}}
\def\RIII{\hbox{${\rm R}_{\rm III}$}}
\subsection{The third case}
We next discuss whether we expect to be able to find a solution
for the third set of parameters~$(a,b,q,x,y)=(1,2,2,8,11)$.
Analogous to before, via linear substitutions and a resultant step,
we should be able to get down to about 6 equations and unknowns,
and we call the resulting system~\RIII.
This already is not the most pleasant computational task,
but only needs be done once (it takes about 15~minutes).
It takes time proportional to~$p^6$ to check all the local solutions,
so we can't take $p$ too much above~20.
The size of the minimal polynomial of a prospective solution
does not matter much due to the quadratic convergence of the Newton method,
but the degree of the field of the solution has a reasonable impact.
We cannot expect to check fields of degree more than 30 or so.
We need $p$ to have a degree~1 factor, but by the Chebotarev density
theorem we can predict that this should happen often enough (even for
a high degree field) so that some prime less than~20 should work.

With these considerations in mind,
we checked the \RIII\ system for local solutions for all primes~$p<20$
and with $p=19$ we found\footnote{
 For both $p=13$ and $p=17$ the local image of this global solution
 was incident with a higher-dimensional solution variety.}
a local solution that lifted to a global \EPZIII\ solution
in the nonic number field given by
$z^9-2z^8-6z^7+8z^6-7z^5+18z^4+44z^3+32z^2+24z+24$,
which has discriminant~$-2^{10}3^7 5^5 11^4$.
For reasons of space, we do not record the solution here.\footnote{
 One model is given modulo~19 by $X(t)=t^8+t^7+6t^6+16t^5+8t^3+4t^2+12$,
 $Q(t)=t^2+3t+13$, $A(t)=16t+15$, $B(t)=17t^2+6t+14$,
 and the interested reader can readily verify that this lifts
 to a $\Q_{19}$-solution with coefficients $x_8,x_7,q_2=1$ and~$y_{10}=0$.}
For \EPZIV\ we were unable to use resultants to reduce
beyond 13 equations and unknowns, and did not even attempt
to find local solutions, even with~$p=5$.
If we had been able to reduce the system down to 10 variables
(as would be hoped from analogy with the above), we could probably
check $p=5$ and maybe~$p=7$.

\section{Concluding comments}
Note that the above four choices of $(a,b,q,x,y)$ are members of
infinite families for which each member has a reasonable possibility
of having infinitely many solutions with~$\rho>10$. Indeed, by taking
\begin{align*}
(a,b,q,x,y)&=(2m,3m,2,10m+2,15m+2)&\rho=10+2/m\\
(a,b,q,x,y)&=(2m,3m+1,2,10m+4,15m+5)&\rho={10m+4\over m+1/3}\\
(a,b,q,x,y)&=(2m+1,3m+1,2,10m+6,15m+8)&\rho={10m+6\over m+1/2}\\
(a,b,q,x,y)&=(2m+1,3m+2,2,10m+8,15m+11)&\rho={10m+8\over m+2/3}
\end{align*}
in each case we have, since $a+b+q+x+y=3x$, the same number of equations
and unknowns, with the value of $\rho=x/\max(a/2,b/3)$ as indicated.
However, we might also suspect that the fields of definition of these
putative solutions become quite large; thus there is no contradiction
with Lang's conjecture, which is only stated for a fixed ground field.

We can also note that with $(a,b,q,x,y)=(2,3,2,10,14)$ we can expect there
to be a nondegenerate $1$-dimensional solution variety~$V$ with~$\rho=10$.
This presumably could be found by a variant of the above methodology,
perhaps by taking specialisations to $0$-dimensional varieties and
finding points on these, and then using this information to reconstruct~$V$.
We have not been able to make this work in practise; although the
specialised system can be reduced to 7 equations and unknowns and we can find
a liftable solution mod~5, it appears that the process of specialisation
increases the degree of the field of the solution beyond our
computational threshold.

\subsection{Performance of computer algebra systems}
In the above computations we used both
\Pari\ \cite{PARI-ref} and~\Magma\ \cite{MAGMA-ref}.
In the end, we were able to do all the relevant computations using
only~\Magma, but this was not apparent at the beginning.
The main difficulty with \Magma\ was dealing with multivariate polynomial
rings, especially as we eliminated variables --- if we did not also
decrease the dimension of the ambient ring, we could experience slowdown.
We also found it to be important to work over the integers rather than
rationals as much as possible,\footnote{
 Except in the small cases where we were able to use the Gr\"obner
 basis machinery; there we want to be working over the rationals rather
 than the integers.}
as else the continual gcd-computations to eliminate denominators
could swamp the calculation.
The availability of multivariate gcd's in \Magma\ frequently allowed us
to reduce the resulting systems by eliminating a common factor.
We found \Magma\ much superior than \Pari\ in searching
for local solutions.\footnote{
 For \RII, \Magma\ took about 5 minutes to find all solutions mod~23,
 and with \RIII\ it took 19 hours to find all solutions mod~19;
 the bulk of the time is actually in computing the determinant of
 the Jacobian matrix to see if the solution lifts, which we could
 ameliorate this partially by (say) not computing the whole Jacobian
 matrix when the first row is zero.}
\Magma\ did quite well in obtaining algebraic numbers
from \hbox{$p$-adic} approximations;
after discussions with the maintainer of \Pari,
we were able to get {\tt algdep} to work
sufficiently well to obtain the above solutions.
The lifting step\footnote{
 We have not made any consideration of the efficacy of a generalisation
 of secant-based methods and/or those of Brent~\cite{brent}.}
was noticeably slower in \Magma\ than in~\Pari,
but as we noted above, the time to do this is not the bottleneck.

\subsection{Acknowledgements}
Thanks are due to Karim Belabas, Nils Bruin, Noam Elkies, and Allan Steel
for comments regarding this work.
The author was partially funded by an NSF VIGRE Postdoctoral Fellowship
at The Pennsylvania State University, the MAGMA Computer Algebra Group
at the University of Sydney, and EPSRC grant GR/T00658/01 during the time
in which this work was done.

\vspace{24pt}
\appendix
\begin{appendix}
\centerline{\textbf{Appendix by Noam D. Elkies} (Harvard University)}
\vspace{8pt}
\centerline
{I. Calculations for the First Case}
\vspace{8pt}
\def\0{^{\phantom0}}

We compute polynomials $X,A,B,Q,Y \in {\bf C}(t)$
of degrees $4,0,1,2,5$, satisfying
\begin{equation}
X^3 + A X + B = Q Y^2.
\label{eq:epz}
\end{equation}
We may normalize $X,Y$\/ to be monic, and translate $t$ so
$Y = t^2 - c$.  Since (\ref{eq:epz}) has degenerate solutions with
$(X,A,B,Y) = (Q(t+b_1)^2, 0, 0, Q(t+b_1)^3)$, we write
\begin{equation}
X = Q ((t+b_1)^2 + 2b_2) + 2 b_3 t + 2 b_4
\label{eq:x}
\end{equation}
for some scalars $b_1,b_2,b_3,b_4$.
Because $AX+B = O(t^5)$ at $t = \infty$,
we have $Y = (X^3/Q)^{1/2} + O(t^{-2})$,
which determines $Y$\/ and imposes two conditions
on $b_1,b_2,b_3,b_4,c$.
Considered as equations in $b_4,c$, these conditions are
simultaneous linear equations, which we solve to obtain
\begin{equation}
b_4 = \frac{b_2^2}{6 b_3} (3 b_3 - 2 b_1 b_2),
\quad
c = \frac
 { (b_3 - b_1 b_2) (3 b_3^2 - 3 b_1 b_2 b_3 + 2 b_2^3)}
 {3 b_2^2 b_3\0}.
\label{eq:b4,c}
\end{equation}
Then $A$ is the $t^4$ coefficient of $Q Y^2 - X^3$; we compute
\begin{equation}
A = \frac{3b_3^2}{b_2^2} (b_3 - b_1 b_2)^2
+ \frac{b_2^2}{3b_3^2} (6 b_1 b_3^3 + 2 b_2^2 b_3^2
  - 6 b_1^2 b_2 b_3^2 - 2 b_1 b_2^3 b_3 + b_1^2 b_2^4
  ).
\label{eq:a}
\end{equation}

The identity (\ref{eq:epz}) then holds if the $t^3$ and $t^2$
coefficients of $X^3 + A X - QY^2$ vanish.
Writing these coefficients in terms of $b_1,b_2,b_3$,
we find that they share a factor $b_3 - b_1 b_2$
that we already encountered in our formula~(\ref{eq:b4,c}) for~$c$.
Namely, the $t^3$ and $t^2$ coefficients are
\begin{equation}
(b_3 - b_1 b_2)
\frac{6 b_3^3 - 6 b_1 b_2 b_3^2 + 6 b_2^3 b_3 - 2 b_1 b_2^4} {3 b_3},
\label{eq:coeff3}
\end{equation}

\begin{equation}
(b_3 - b_1 b_2)
\frac
{18 b_3^5 + (15 b_2^3 - 18 b_1^2 b_2^2) b_3^3 + 15 b_1\0 b_2^4 b_3^2
 + (2 b_2^6 - 6 b_1^2 b_2^5) b_3\0 - 2 b_1\0 b_2^7}
{9 b_2\0 b_3^2}.
\label{eq:coeff2}
\end{equation}
If $b_3 = b_1 b_2$ then $c=0$ and $b_4 = b_2^2/6$,
and we calculate $A = -b_2^4/3$ and $B = 2 b_6^2 / 27$.
But this makes $X^3+AX+B = (X+2(b_2^2/3))(X-(b_2^2/3))^2$,
so our elliptic curve degenerates to a rational curve with a node
(or a cusp if $b_2$ vanishes too).

Therefore the numerators of the fractions
in (\ref{eq:coeff3},\ref{eq:coeff2}) must vanish.
The first of these yields a linear equation in~$b_1$,
which we solve to obtain
\begin{equation}
b_1 = \frac{3 b_3 (b_3^2 + b_2^3)} {b_2 (3 b_3^2 + b_2^3)}.
\label{eq:b1}
\end{equation}
Substituting this into (\ref{eq:coeff2}) yields
$2 b_2^6 b_3\0 (3 b_3^2 - 2 b_2^3) / (3 b_3^2 + b_2^3)$.
We conclude that $3b_3^2 = 2 b_2^3$.

All nonzero solutions of $3b_3^2 = 2 b_2^3$ are equivalent
under scaling.  We choose $(b_2,b_3)=(6,12)$ and work our way back.
We find $b_1=10/3$, and then $c=-8/9$, $b_4=-2$, and finally
$A=528$ and $B=128(12t+31)$.  To optimize the constants
in the resulting family of large integral points on elliptic curves,
we replace $t$ by $6t-(10/3)$ and renormalize to obtain at last
\[ A=33, \; B=-18(8t-1), \; Q=9t^2-10t+3 , \]
% arXiv formatting is different as usual
\vskip0pt minus 6pt
\begin{equation}
X = 3 (108 t^4 - 120 t^3 + 72 t^2 - 28 t + 5)
\end{equation}
\[ Y = 36 (54 t^5 - 60 t^4 + 45 t^3 - 21 t^2 + 6 t - 1) .  \]

To complete the proof that there are no other solutions,
we must also consider the possibility that the denominator
of (\ref{eq:b4,c}) vanishes, which is to say $b_2=0$ or $b_3=0$.
If $b_2=0$ then the $t^6$ and $t^5$ coefficients of $X^3-QY^2$
reduce to $3b_3^2$ and $6 b_3 (b_1 b_3 + b_4)$.
Thus we also have $b_3=0$, and then
$A=-3b_4^2$ and $X^3+AX-QY^2=2b_4^3$,
so the condition on the $t^3$ and $t^2$ coefficients
holds automatically for any choice of $b_4$ and~$c$.
But this makes $X^3+AX+B = (X-2b_4)(X+b_4)^2$,
so again we have a degenerate elliptic curve.
If $b_3=0$ but $b_2 \neq 0$ we obtain
$b_1=0$ and $6 b_4 = b_2 (3c - b_2)$.
Then $A = -b_2^2 (9 c_2\0 + 4 b_2^2) / 12$, and $X^3+AX-QY^2$
has $t^3$~coefficient zero but $t^2$~coefficient $b_2^3 c^2/2$.
Since we assume $b_2\neq 0$, we conclude $c=0$,
leaving $A=-b_2^4/3$ and $B = 2 b_2^6/27$,
for the same degenerate elliptic curve as above.
\goodbreak
\centerline{II. Letter from Noam D. Elkies to Don Zagier (1988)}
\vspace{8pt}
\noindent{Dear Prof.\ Zagier,}
\vspace{8pt}

I have read with considerable pleasure your note on ``Large integral 
points on elliptic curves'', which Prof.~Gross showed me in response
to a question.  In the second part of that note you you define a
``measure of impressiveness'', $\rho$, \,of a large integral point
$(x,y)$ on the elliptic curve $x^3+ax+b=y^2$ by 
\[ \rho =\log (x)/\log (\max (|a|^{\frac{1}{2}},|b|^{\frac{1}{3}}))\] 
and exhibit several infinite families of such points for which 
$\rho =9+O(\frac{1}{\log x})$.  You conjectured, though, that $\rho $
could be as large as 10, so I searched for an infinite family
confirming this.  What I found was an infinite family of Pell type
for which $\rho =12-O(\frac{1}{\log x})$.  The implied constant is
quite large---bigger than 200---so $\rho $ approaches 12 very slowly,
remaining below $5\frac{1}{2}$ for $x$ in the range $[1,10^8]$ of
Odlyzko's computation, and first exceeding 10 and 11 for $x$ of 51
and 107 digits respectively.

In your note you give a probabilistic heuristic suggesting that
$\rho $ should never significantly exceed 10.  But a na\"{\i}ve
counting of parameters and constraints for a Pell-type family
\begin{equation}\label{Z1}
X^{3}(t)+A(t)X(t)+B(t)=Q(t)Y^{2}(t)       
\end{equation}
(in which $A$, $B$ are polynomials of low degree, $Q$ is a quadratic
polynomial in $t$, and $X$, $Y$ are polynomials of large degree)
suggests that \eqref{Z1} should have several solutions with $\rho 
\rightarrow 12$, most simply with $A$ constant, $B$ linear, $X$
quartic and $Y$ quintic.  Actually finding such a solution required
a longish MACSYMA session to solve four nonlinear equations in four
variables, which surprisingly have a unique nontrivial solution,
(necessarily) defined over ${\bf Q}$: up to rescaling $t$ and 
the polynomials $A$, $B$, $Q$, $X$, $Y$, the only solution to
\eqref{Z1} is
\[ A=33, \: B=-18(8t-1), \: Q=9t^{2}-10t+3 , \]
\begin{equation}\label{Z2}
X=324t^{4}-360t^{3}+216t^{2}-84t+15 ,
\end{equation}
\[ Y=36(54t^{5}-60t^{4}+45t^{3}-21t^{2}+6t-1) . \]
As it stands, \eqref{Z2} seems of little use because $Q$ is never a 
square for $t\in {\bf Z}$.  However, we may rescale \eqref{Z2} by replacing
$(A, B, X)$ by $(4A=132, 8B, 2X)$, which yields an integral
point provided $2Q$ is a square.   That Pell-type condition is
satisfied by $t=1$ and thus by infinitely many $t$, yielding an
infinite family of solutions $(b, x, y)$ to
$x^{3}+132x+b=y^{2}$
with $x\sim 2^{-25}3^{-4}b^{4}$.  The small factor
$2^{-25}3^{-4}\doteq 3.68\cdot 10^{-10}$ means that, although $\rho $
eventually approaches 12, the first few admissible values of $t$
yield only mediocre $\rho $: the second such value, $t=15$, when
$b=-17424$ and $x=35334750$ (the largest such $x$ to fall within
the bounds of Odlyzko's search), produces only {$\rho \doteq 5.34$}
and was probably ignored; only the ninth value \ $t=812111750209$
produces $\rho >10$, and only the eighteenth, 
$t=-48926085100653611109021839$, reaches $\rho >11$. 

Some final remarks:  Prof.~Lang tells me that Vojta's conjectures
imply the $\rho \leq 10+\epsilon$ conjecture {\em except possibly
for a finite number of exceptional families} such as those
obtained by rescaling \eqref{Z2}.  Vojta proves this implication in a yet
unpublished paper, but leaves open the existence of exceptional
families.  It's interesting to compare this situation with the
similar conjecture of Hall concerning $|x^{3}-y^{2}|$, where the  
best infinite families known come from the identity
\begin{equation}\label{Z3}
(t^{2}+10t+5)^{3}-(t^{2}+22t+125)(t^{2}+4t-1)^2=1728t 
\end{equation}
(Exer.\   9.10 in Silverman's {\em The Arithmetic of Elliptic Curves},
attributed to Danilov, {\em Math.\ Notes Acad.\ Sci.\ USSR} {\bf 32}
(1982), 617--8), which yields Pell-type solutions with $\rho $
tending this time to the ``correct'' value of 6.  There is a 
natural reason (which Danilov does not mention in his article) 
for \eqref{Z3} to be defined over ${\bf Q}$: the fifth modular
curve $(j(z),j(5z))$ is rationally parametrized by
\[
j(z)=f(t)=\frac{(t^{2}+10t+5)^{3}}{t}, \: j(5z)=f(\frac{1}{t}),
\]
and $f(t)$ is a sixth-degree rational function with a fifth-order
pole at infinity (a cusp), two third-order zeros (CM by 
$\frac{1}{2}(1+\sqrt{-3})$) and two second-order values of 1728
(CM by $\sqrt{-1}=i$; the appearance of $z=\frac{1}{5}(i\pm 2)$
when $j(z)=j(5z)=1728$ splits the other two inverse images of
1728 under~$f$)---hence~\eqref{Z3}.  I have no similar rationale
for~\eqref{Z2}, nor for why it gives ``too large'' a value of~$\rho $.
\hfil\break\break
\vspace{2pt}\noindent
\hspace{3.5in}Sincerely,\hfil\break
\vspace{2pt}\noindent
\hspace{3.5in} (signed)\hfil\break
\vspace{0pt}\noindent
\hspace{3.5in}Noam D. Elkies
\vspace{24pt}
\end{appendix}


\begin{thebibliography}{99}
\bibitem{macsyma-ref}
\textsc{Macsyma}, a sophisticated computer algebra system.
See \verb+maxima.sourceforge.net+ for history and current version
of its descendants.
\end{thebibliography}
\end{document}